\documentclass[11pt]{article}
\usepackage[pctex32]{graphics}
\usepackage{amssymb}
\usepackage{latexsym}
\usepackage{epsfig} 
\usepackage{multirow} 
\usepackage{pstricks}
\usepackage{graphicx}
\usepackage{float}
\usepackage{epstopdf}
\usepackage{mathrsfs}
\usepackage{subcaption}
\usepackage{soul}

\usepackage{amsmath}
\usepackage{amssymb}
\usepackage{relsize}
\usepackage{multirow}
\usepackage{color}
\usepackage{graphicx}
\usepackage{multicol} 

\usepackage{times}
\usepackage{epstopdf}
\newcommand{\comm}[1]{}
\newtheorem{definition}{Definition}[section]
\newtheorem{theorem}[definition]{Theorem}
\newtheorem{lemma}[definition]{Lemma}

\usepackage{fancyhdr}
\usepackage{setspace}
\usepackage{pstricks} 
\newcommand{\R}{{\mathbb R}}

\newcommand{\mM}{{\mathsf M}}

\newcommand{\mS}{{\mathsf S}}

\newcommand{\mP}{{\mathsf P}}
\newcommand{\mR}{{\mathsf R}}
\newcommand{\mT}{{\mathsf T}}
\newcommand{\mC}{{\mathsf C}}
\newcommand{\mI}{{\mathsf I}}
\newcommand{\mA}{{\mathsf A}}
\newcommand{\mV}{{\mathsf V}}

\newcommand{\mB}{{\mathsf B}}
\newcommand{\mL}{{\mathsf L}}
\newcommand{\mQ}{{\mathsf Q}}

\newcommand{\mD}{{\mathsf D}}

\newcommand{\mU}{{\mathsf U}}
\newcommand{\mG}{{\mathsf G}}

\newcommand{\mSigma}{{\mathsf \Sigma}}

\newcommand{\mGamma}{{\mathsf \Gamma}}

\begin{document} 

\title{Subspace Splitting Fast Sampling from Gaussian Posterior Distributions of Linear Inverse Problems}
\author{Daniela Calvetti\and Erkki Somersalo}
\date{ Case Western Reserve University \\
Department of Mathematics, Applied Mathematics, and Statistics}

\maketitle
\begin{abstract}
It is well-known that the posterior density of linear inverse problems with Gaussian prior and Gaussian likelihood is also Gaussian, hence completely described by its covariance and expectation. Sampling from a Gaussian posterior may be important in the analysis of various non-Gaussian inverse problems in which a estimates from a Gaussian posterior distribution constitute an intermediate stage in a Bayesian workflow. Sampling from a Gaussian distribution is straightforward if the Cholesky factorization of the covariance matrix or its inverse is available, however when the unknown is high dimensional, the computation of the posterior covariance maybe unfeasible. If the linear inverse problem is underdetermined, it is possible to exploit the orthogonality of the fundamental subspaces associated with the coefficient matrix together with the idea behind the Randomize-Then-Optimize approach to design a low complexity posterior sampler that does not require the posterior covariance to be formed. 
The performance of the proposed sampler is illustrated with a few computed examples, including non-Gaussian problems with non-linear forward model, and hierarchical models comprising a conditionally Gaussian submodel.
\end{abstract}

\section{Introduction}

The popularity of Gaussian distributions in applications arises from the wealth of their well-known properties, including their concise and complete parametric description in terms of expectation and covariance matrix. In the context of inverse problems within the Bayesian framework, it is well known that the posterior density associated with linear inverse problems with Gaussian priors and additive Gaussian errors is a Gaussian distribution, with  mean and covariance expressible in terms of the mean and covariance of the prior and noise, and of the forward operator matrix and the observed data constituting the right hand side of the linear system.  

In some applications, despite the existence of the closed form solution, it is necessary to augment the analytic description of the density with a sample of realizations, e.g., when the covariance matrix is too large to be computed or the forward model is given in a matrix-free form. Other instances when sampling is desirable is to obtain marginal uncertainty estimates, or to draw scatter plots of few components or derived quantities. Sampling from a Gaussian density is necessary also in inverse problems that include an intermediate linear Gaussian inverse problem in a workflow of Bayesian analysis, e.g., constituting a proposal for and MCMC algorithm \cite{cotter2013mcmc}, or when uncertainty quantification of the absence/presence of an anomaly such as a malignancy in an image is based on statistical hypothesis testing \cite{rambojun2024uncertainty}.
Estimating a non-linear function of a Gaussian random variable with uncertainty analysis is an example in which samples from the Gaussian density are useful, e.g., in the estimation of quadratic heat absorption, or thermal dose, in ultrasound surgery \cite{malinen2004simulation}. 
In uncertainty quantification of non-Gaussian hierarchical models, effective sampling from a conditionally Gaussian density may be necessary \cite{Calvetti2025dictionary}.
Finally, the sampling from Gaussian posteriors constitutes an essential part of the analysis step of Ensemble Kalman Filtering (EnKF) algorithm and Ensemble Kalman inversion (EKI) \cite{iglesias2021adaptive}.

In principle, sampling from Gaussian distributions is straightforward if a symmetric factorization such as Cholesky factorization of either the covariance or the precision matrix is available and can be easily accessed, but may become prohibitively expensive when the dimension of the random variable is very large and the computing environment lacks the computational resources required for formation and factorization of the covariance matrix. In the literature, there are multiple papers suggesting how to overcome the curse of dimensionality in sampling using different tools, ranging from iterative numerical linear algebra to Markov chain Monte Carlo (MCMC), characteristic of the proposing community. A recent review and comparison of these methods can be found in \cite{vono2022high}. 

Due to their central role in computations, symmetric factorizations of symmetric positive definite (SPD) matrices have been studied extensively in numerical linear algebra. It is well-known that the computational complexity and memory allocation required for a Cholesky factorization decrease sharply for sparse and banded matrices, therefore it has been suggested in \cite{rue2001fast} to rearrange the order of the unknowns so as to cluster the nonzero entries of the covariance or precision matrix around the main diagonal. Since the symmetric tridiagonal factorization of a symmetric positive definite matrix can be obtained from the analogous factorization of its inverse, the natural choice when trying to keep the computational complexity low is to compute the Cholesky factorization of the matrix that is sparser between the covariance and its inverse, known as the precision.

In the context of interest in this paper, the target density arises from a linear observation model of the form 
\begin{equation}\label{linsys}
b = \mA x + \varepsilon, \; x\in\R^n, \; b\in\R^m,
\end{equation}
where $x$ is the unknown parameter of interest, $\mA$ is a known real matrix, and $\varepsilon$ is additive observation noise. Following the Bayesian paradigm, $x$ and $\varepsilon$ are modeled as random variables, which we assume to be mutually independent and normally distributed, i.e.,  
\begin{equation}\label{standard normal}
 x\sim \pi_x(x) = {\mathcal N}(x\mid x_0,\mGamma), \quad \varepsilon \sim{\mathcal N}(0,\mSigma).
\end{equation}
It follows from Bayes' formula \cite{calvetti2023bayesian} that under these assumptions the posterior distribution, given by
\begin{eqnarray*}
 \pi_{x\mid b}(x\mid b) &\propto& \pi_x(x) \pi_{b\mid x}(b\mid x) \\
 &\propto& {\rm exp}\left( - \frac 12 (x-x_0)^\mT \mGamma^{-1} (x - x_0) - \frac 12 (b - \mA x)^\mT\mSigma^{-1}(x - \mA x)\right),
\end{eqnarray*}
 is also a Gaussian distribution, 
\begin{equation}\label{post x}
 \pi_{x\mid b}(x\mid b) = {\mathcal N}(x\mid \mu,\mC),
\end{equation}
with posterior covariance and expectation of the form  
\begin{equation}\label{posterior_cov}
 \mC = \big(\mA^\mT \mSigma^{-1} \mA +\mGamma^{-1}\big)^{-1}, \quad \mu = \mC (\mA^\mT \mSigma^{-1} b +\mGamma^{-1}x_0).
\end{equation}
To characterize the posterior distribution by means of a single realization, it is common to compute the Maximum a Posteriori (MAP) estimate, which in the Gaussian case and with standard parametrization of the model coincides with the posterior mean and can be found by minimizing the negative logarithm of the posterior, 
\begin{equation}\label{MAP}
x_{\rm MAP} = {\rm argmin}\big\{ (b-\mA x)^\mT \mSigma^{-1} (b - \mA x) + (x - x_0)^\mT \mGamma^{-1} (x - x_0)\big\}.
\end{equation}
Furthermore, if the Cholesky factorization, or any other symmetric factorization, of the posterior covariance (\ref{posterior_cov}) is available,
\[
\mC = \mR^\mT \mR, 
\]
independent sampling from the posterior is rather straightforward. In fact, if $w^{(j)} \sim {\mathcal N}(0, \mI_n)$ then 
\[
x^{(j)} = \mR^{-\mT}w^{(j)} + \mu 
\]
is a realization from the posterior.
Alternatively, if we have access to a symmetric factorization of the posterior precision matrix $\mP = \mC^{-1}$,
\begin{equation}\label{posterior_prec}
\mP = \mA^\mT\mSigma^{-1}\mA +\mGamma^{-1} = \mG^\mT \mG, 
\end{equation}
recalling that, in this case,
\[
\mC = \mG^{-1} \mG^{-\mT},
\]
we have that if $w^{(j)} \sim {\mathcal N}(0, \mI_n)$ then 
\[
x^{(j)} = \mG w^{(j)} + \mu 
\]
is a realization from the posterior.
This direct approach to sampling may not be possible if the direct computation of the posterior covariance from (\ref{posterior_cov}) or the posterior precision from (\ref{posterior_prec}) are prohibitively expensive.

The matrix $\mA^\mT \mSigma^{-1} \mA $ is the Hessian of the objective function minimized when fitting the data, and referred to as the Hessian of the data misfit. It has been observed that in many applications the Hessian portion of the posterior covariance is of low rank. In \cite{flath2011fast} the authors propose to take advantage of the low dimensionality and possible structure of the Hessian to reduce the computational complexity of drawing from the Gaussian posteriors of large scale linear inverse problems with Gaussian likelihood and prior. 

In the case where the dimensionality of the data space $\R^m$ is smaller than that of the unknowns, it may advantageous to consider the adjoint problem and perform to computationally demanding operations in data space rather than in the input space of the unknown: see, e.g., \cite{calvetti2000regularizing,sanz2023inverse}. In this case, expressing the posterior covariance and mean in the alternative form 
\begin{equation}\label{post_cov2}
\mC = \mGamma - \mGamma \mA^\mT\big( \mA \mGamma \mA^\mT + \mSigma \big)^{-1} \mA \mGamma, \quad\mu = x_0 +
\mGamma \mA^\mT\big( \mA \mGamma \mA^\mT + \mSigma \big)^{-1}b,
\end{equation}
may result in a substantial reduction in the computational complexity for its calculation, because the matrix to be inverted is of the dimension of the data, as suggested also in \cite{kaipio2006statistical,sanz2023inverse}. 
In general, the relative dimensions of the data and the unknown can be used to select how to express the posterior covariance in terms of the prior and noise covariance, preferring 
(\ref{post_cov2}) in the underdetermined case, and (\ref{posterior_cov}) in the overdetermined case. If the dimensionality of the unknown is high, the posterior covariance is large regardless of how it is computed. In this article, the main focus is answering the following question: {\em If the forward model is strongly underdetermined, $m\ll n$, how can we make use of the low dimensionality of the data space in random sampling?} 

Central to the discussion in this paper is the idea of using optimization together with a perturbation of the data and the prior mean, proposed in \cite{chen2012ensemble} in the context of Ensemble Kalman filtering. In  \cite{bardsley2014randomize}, it is shown that the samples produced with this scheme, known as Randomize-Then-Optimize (RTO), are distributed according to the posterior distribution if the observation model (\ref{linsys}) is linear. For non-linear models, the perturbation argument is used to generate proposals for the Metropolis-Hastings algorithm. The details of the RTO algorithms and its properties will be reviewed further below. In this work, using an appropriate subspace splitting, we propose a two-phase optimization-based sampling scheme that requires calculations only in the low-dimensional data space.

The rest of the paper is organized as follows. In section 2 we review the key ideas of the RTO where a sample from the Gaussian posterior associated with a linear problem with Gaussian prior and Gaussian likelihood is constructed by perturbing the data and the prior mean. To simplify the notation, in this section we assume that both the prior and the likelihood are zero mean white Gaussian. 
We then proceed to illustrate how, using the properties of the fundamental subspaces associated with the coefficient matrix, it is possible to design a computationally efficient algorithm for the RTO sampler suitable for matrix-free problems. In the same section we derive a formulation of the algorithm based on the adjoint of the linear operator which significantly reduces the computational complexity for strongly underdetermined problems, $m\ll n$.
In section 3 we show how the results of section 2 extend to the general Gaussian setting for linear problems, where the covariances of the likelihood and of the prior are different from the identity, and to the case where the white Gaussian prior is for a random variable $z$ which is a linear function of the unknown of interest $x$, 
\[
z = \mL x, \; \mL \in \R^{p \times n}, \; p\geq n,
\]
with $\mL$ of full rank $n$. Section 4 show how the computations of the sampling algorithm presented in section 2 can be organized in the underdetermined case to reduce the complexity of the algorithm.
In section 5 we outline how Krylov subspace iterative linear solvers can be employed when the matrix $\mA$ is not necessarily explicitly available but it is possible to compute products of arbitrary vectors with $\mA$ and $\mA^\mT$. 
In section 6 the performance of the proposed algorithms is illustrated with a couple of computed examples, including an extension of the proposed method to non-linear, hence non-Gaussian, inverse problems, where the sampling from the Gaussian distribution serves as an effective way of generating proposals for the preconditioned Crank-Nicholson (pCR) sampling scheme. Furthermore, sampling from a non-Gaussian posterior with a conditionally Gaussian prior model is discussed in the light of a dictionary learning example. 

\section{Splitting-based Randomize-Then-Optimze sampling of a Gaussian posterior}

We begin by reviewing the well-known random perturbation argument that connects tightly optimization and sampling in the Gaussian case. 
In the following, we use the notation $\|z\|_\mM^2 = z^\mT\mM^{-1} z$, where $\mM$ is a given SPD matrix.
\begin{theorem}\label{th:get xj}
Assume that $\mA\in\R^{m\times n}$, and $b\in\R^m$ is an indirect noisy observation (\ref{linsys}) of $x\in\R^n$, where
\begin{equation}\label{gaussian model}
 x \sim{\mathcal N}(x_0,\mGamma), \quad \varepsilon\sim{\mathcal N}(0,\mSigma),
\end{equation}
and $x$ and $\varepsilon$ are mutually independent. We define
\begin{equation}\label{get xj}
 x^j = {\rm argmin}\big\{ \|b^j - \mA x\|_\mSigma^2 + \|x - x_0^j\|_\mGamma^2\big\}, \quad j=1,2,\ldots
\end{equation}
where $b^j = b+\eta^j$ and $x_0^j = x_0 + \nu^j$, the random perturbations being drawn from the densities
\[
 \eta^j\sim{\mathcal N}(0,\mSigma), \quad \nu^j\sim{\mathcal N}(0,\mGamma).
\]
Then the points $x^j$ are independent samples from the posterior density $\pi_{x\mid b}(x\mid b)$.
\end{theorem}

{\bf Proof:} It follows from assumptions of the theorem that the posterior mean estimate which coincides with the MAP estimate for the Gaussian posterior (\ref{posterior_cov}), is of the form
\[
 x^j =  \mC\big(\mA^\mT\mSigma^{-1} b^j +\mGamma^{-1}x_0^j\big).
\]
Since the random perturbations have a zero mean, we see that
\[
\overline x =
 {\mathbb E}(x^j) = \mC\big(\mA^\mT\mSigma^{-1} b +\mGamma^{-1}x_0\big) = \mu,
\]
which is the posterior mean. Moreover,
\begin{equation}\label{post uncertainty}
 x^j - \overline x  = \mC\mA^\mT\mSigma^{-1} (\eta^j + \mGamma^{-1}\nu^j),
\end{equation}
implying that the covariance matrix of the sample is
\begin{eqnarray*}
{\mathbb E}\big(( x^j - \overline x)(x^j - \overline x)^\mT\big) &=&
  \mC\big\{\mA^\mT\mSigma^{-1}
{\mathbb E}\big(\eta^j(\eta^j)^\mT\big)\mSigma^{-1}\mA +\mGamma^{-1}{\mathbb E}\big(\nu^j(\nu^j)^\mT\big)\mGamma^{-1}
\big\} \mC  \\
&=& \mC\big\{\mA^\mT\mSigma^{-1} \mA+\mGamma^{-1}\big\}\mC =\mC,
\end{eqnarray*}
proving that the sample is distributed according to the posterior. $\Box$

One important observation is that the posterior uncertainty, given by formula (\ref{post uncertainty}) can be split in two parts, one arising  from the prior, the other from the likelihood. The second point is that in the Gaussian case posterior sampling can be formulated as an optimization problem (\ref{get xj}). The latter observation is key for the remainder of the paper, which addresses how linear algebra can be used to make this optimization computationally effective.

In the rest of this section, we assume that the covariances of the prior and likelihood (\ref{standard normal}) are identity matrices, i.e.,
\[
\mGamma = \mI_n, \quad \mSigma = \mI_m, 
\]
and that the prior mean vanishes, i.e., $x_0=0$. Later, the modifications necessary for the general case will be discussed in detail. With these simplifying assumptions,  the posterior covariance and expectation can be expressed as
\[
\mC = \big(\mA^\mT \mA+ \mI_n\big)^{-1}, \quad \mu = \mC\mA^\mT b
\]
and the posterior sampling of Theorem~\ref{th:get xj} can be reduced to solving the linear least squares problem
\begin{equation}\label{x_min}
 x = {\rm argmin}\big\{ \|b + \eta - \mA x\|^2 + \|x - \nu\|^2\big\}, \quad \eta\sim{\mathcal N}(0,\mI_m), \quad \nu \sim{\mathcal N}(0,\mI_n).
\end{equation}
Here and in the following discussion, to simplify the notation we omit the superscript $j$.
By writing
\[
 \|b + \eta - \mA x\|^2 + \|x - \nu\|^2 =\left\|
 \left[\begin{array}{c} \mA \\ \mI_n\end{array}\right] x - \left[\begin{array}{c} b + \eta \\ \nu\end{array}\right]
 \right\|^2,
\]
it follows that the minimizer of (\ref{x_min}) is the solution of the normal equations
\begin{equation}\label{ATA perturbed}
 \big(\mA^\mT\mA + \mI_n\big)x = \mA^\mT(b + \eta) + \nu.
\end{equation}
In general, if the problem is high dimensional, the solution of the normal equations can be computed , e.g., by using Krylov subspace based iterative solvers, as will be discussed later in this paper. It is worth noticing that if $m\ll n$, the data space has much lower dimension than the parameter space and significant computational advantages can be obtained by arriving at the solution of (\ref{ATA perturbed}) through the solution of an adjoint problem in data space.

\subsection{Splitting-based sampling and adjoint problem}

In this section, we assume that $m<n$. We begin by stating the following result.

\begin{lemma}\label{Lemmma22}
Let $x\in\R^n$ be the unique solution of the normal equations
\begin{equation}\label{ATA}
 (\mA^\mT\mA + \mI_n) x = \mA^\mT b,
\end{equation}
and let $z\in\R^m$ be the unique solution of the adjoint problem
\begin{equation}\label{AAT}
 (\mA\mA^\mT + \mI_m) z = b.
\end{equation}
Then, $x = \mA^\mT z$.    
\end{lemma}

{\bf Proof:} Multiplying equation (\ref{AAT}) from the left by $\mA^\mT$ gives
\[
 \mA^\mT b = \mA^\mT\big(\mA\mA^\mT + \mI_m\big) z 
 = \big(\mA^\mT\mA + \mI_n\big)\mA^\mT z,
\]
and the claim follows from the uniqueness of the solution of (\ref{ATA}). $\Box$

It follows from the Lemma above that when $m<n$, the complexity of the solution of problem (\ref{ATA}) can always be reduced by solving the lower dimensional problem (\ref{AAT}). This observation can be used to speed up the computation of large scale  Tikhonov-regularized solutions of underdetermined problems, see, \cite{calvetti2024distributed}.
Our aim in this paper is to extend these complexity reduction ideas to sampling schemes. A complication that needs to be addressed is the fact that
in the case of sampling, the right hand side of the equation (\ref{ATA perturbed}) is not in the range of the operator $\mA^\mT$, hence requiring an appropriate splitting of the prior noise term $\nu$.

The relations between the fundamental subspaces of a matrix implies that the term $\mA^\mT(b + \eta)$ on the right hand side of (\ref{ATA perturbed}) is orthogonal to the null space of $\mA$. Next we decompose the perturbation the term $\nu$ as the sum of two orthogonal components. 

\begin{lemma}\label{lemma:reduce}
The term $\nu$ in (\ref{ATA perturbed}) can be expressed as 
\begin{equation}\label{partition nu}
 \nu = \mA^\mT \delta  + h 
\end{equation}
where $h \in {\mathcal N}(\mA)$ and $\delta$ is the least squares solution of the overdetermined linear system 
\begin{equation}\label{LSQ delta}
 \mA^\mT \delta = \nu,
\end{equation} 
that is,
\[
 \delta = {\rm argmin} \big\{\|\mA^\mT\delta-\nu\|^2\big\}.
\]
\end{lemma} 

{\bf Proof.}
The existence of the partitioning (\ref{partition nu}) follows directly from the observation that $\R^n = {\mathcal R}(\mA^\mT) \oplus {\mathcal N}(\mA).$

Multiplying both sides of (\ref{partition nu}) by $\mA$ from the left shows that $\delta$ must satisfy 
\[
 \mA \mA^\mT \delta = \mA \nu,
\]
which are the the $m \times m$ normal equations associated with the overdetermined linear system
(\ref{LSQ delta}),  thus completing the proof. $\Box$

If the decomposition (\ref{partition nu}) is available,
we express $x$ in (\ref{ATA perturbed}) as 
\[
 x = x_1 + x_2, \quad x_1\in {\mathcal R}(\mA^\mT), \quad x_2\in {\mathcal N}(\mA),
\]
and writing the right hand side of (\ref{ATA perturbed}) in terms of (\ref{partition nu}),
we have
\begin{eqnarray*}
  \big(\mA^\mT\mA + \mI_n\big)(x_1 + x_2) &=&
  \big(\mA^\mT\mA + \mI_n\big)x_1 + x_2  \\
  &=&
  \mA^\mT(b + \eta +\delta) + h.
\end{eqnarray*}
It follows from the orthogonality of the fundamental subspaces that $x_1$ and $x_2$ must satisfy
\begin{equation}\label{normal x1}
 (\mA^\mT \mA + \mI_n)  x_1  = \mA^\mT (b +  \eta +  \delta), \quad x_2 = h, 
\end{equation}
hence Lemma~\ref{lemma:reduce} give us the following result.

\begin{theorem}\label{th:split}
Given the decomposition (\ref{partition nu}), the unique solution of the linear system (\ref{ATA perturbed}) can be written as
\[
 x = \mA^\mT z + h,
\]
where $z\in\R^m$ is the unique solution of the adjoint equation
\begin{equation}\label{get z}
 \big(\mA\mA^\mT + \mI_m\big) z = b+\eta+\delta.
\end{equation}  
\end{theorem}

We summarize the previous results in the following theorem outlining how the computations for a Gaussian posterior sampling scheme can be organized.

\begin{theorem}
Consider the problem of generating an independent sample $\{x^1,x^2,\ldots,x^K\}$ drawn from the posterior distribution of the whitened linear model,
\[
\pi_{x\mid b}(x\mid b) = {\mathcal N}(x\mid \mu,\mC), \quad \mC = \left(\mA^\mT\mA + \mI_n\right)^{-1}, \quad \mu = \mC \mA^\mT b,
\]
where $\mA\in\R^{m\times n}$ has full rank.
\begin{enumerate}
\item For each $j$, draw independent perturbations
\[
 \nu^j \sim{\mathcal N}(0,\mI_n), \quad \eta^j \sim{\mathcal N}(0,\mI_m).
\]
\item If $n\leq m$, solve the normal equations
\begin{equation}\label{normal n}
 \big(\mA^\mT\mA + \mI_n\big)x^j = \mA^\mT(b + \eta^j) + \nu^j.
\end{equation}
\item If $m<n$,
\begin{enumerate}
    \item Solve the normal equations
    \begin{equation}\label{get delta}
 \mA \mA^\mT \delta^j = \mA \nu^j,
\end{equation}
and define $h^j = \nu^j - \mA^\mT\delta^j$.
\item Solve the normal equations
\begin{equation}\label{adjoint m}
 \big(\mA\mA^\mT + \mI_m\big) z^j = b+\eta^j+\delta^j.
\end{equation}
\item Set
\[
 x^j = \mA^\mT z^j + h^j.
\]
\end{enumerate}
\end{enumerate}
\end{theorem}

Before discussing the numerical algorithms to solve the various least squares problems, we spend a few words on the well-posedness of the step 3 (a) above, and on the whitening of the linear forward model.

The assumption that the matrix $\mA$ has full rank does not prevent the non-zero singular values from ranging over orders of magnitudes, implying that step  3 (a) may be of concern because of of its dependency on the square of the condition number of $\mA$. Fortunately, the intermediate variable $\delta^j$ is used only to compute an orthogonal projection of $\nu^j$ onto the null space of $\mA$. Introducing the reduced singular value decomposition (SVD) of the matrix $\mA$, 
\[
 \mA = \sum_{k=1}^m \sigma_j u_j v_j^\mT,
\] 
we observe that
\[
 \delta^j = \sum_{k=1}^m \frac{(v_k^\mT\nu^j)}{\sigma_k} u_k,
\]
hence, the contribution to $h^j$ is
\[
 \mA^\mT\delta^j = \sum_{k=1}^m \big(v_k^\mT\nu^j\big) v_k.
\]
Similarly, the contribution of $\delta^j$ to $x^j$ through step 3 (b) is given by 
\[
 \mA^\mT\big(\mA\mA^\mT + \mI_m\big)^{-1}\delta_j = 
 \sum_{k=1}^m \frac{\big(v_k^\mT\nu^j\big)}{1+\sigma_k^2} v_k.
\]
Therefore, while the normal equations (\ref{get delta}) may be ill-conditioned and as such may require standard regularization, the effect of the ill-conditioning on $x^j$ is mitigated by the multiplications by $\mA^\mT$.
The numerical solution of these equations will be discussed in detail later.

\subsection{Matrix-free sampling }

In some applications the forward operator is given in matrix-free form, in the sense that the matrix $\mA$ is not explicitly available, but it is possible to compute the product of $\mA$ and $\mA^\mT$ with arbitrary vectors. Alternatively, the matrix $\mA$ may be large, and solving linear systems involving it using direct methods may not be feasible. It is rather straightforward to implement the proposed sampling scheme so that it can be used with matrix-free forward models. In this section we provide an overview of how the calculations can be carried out without requiring direct access to the matrix.

Given a matrix  $\mA\in\R^{m\times n}$, a vector $b\in \R^m$, and a parameter $\mu\geq 0$, consider the following two linear systems:
\begin{enumerate}
\item The normal equations, where we seek a vector $x\in\R^n$ satisfying
\begin{equation}\label{normal}
\big(\mA^\mT\mA + \mu\, \mI_n\big) x = \mA^\mT b,
\end{equation}
\item The adjoint equations, where we seeks a vector $z\in\R^m$ that solve
\begin{equation}\label{adjoint}
\big(\mA\mA^\mT + \mu\, \mI_m\big) z = b.
\end{equation}
\end{enumerate}
It follows from Lemma~\ref{Lemmma22} that when $\mu>0$ each problem has a unique solution because the coefficient matrices are invertible, and solutions $x$ and $z$ satisfy the linear system
\[
 x = \mA^\mT z.
\] 
The Lanczos bidiagonalization process, originally proposed by Lanczos \cite{lanczos1952solution}, also known as the Golub-Kahan algorithm \cite{golub2009matrices}, provides a unified framework for dealing with both types of equations regardless of the relative sizes of $m$ and $n$. Below, we review the process and indicate how it can be used to find an approximate solution for both types of equations. 

The Lanczos process generates iteratively a sequence of matrix triplets $\big\{\big(\mU_\ell,\mB_\ell,\mV_\ell\big)\mid \ell = 1,2,\ldots\big\}$ with the following properties:

\begin{enumerate}
\item $\beta_1 u^{(1)} = b$,
\[
 \mU_\ell = \left[\begin{array}{ccc} u^{(1)} & \cdots & u^{(\ell)}\end{array}\right], \quad \mbox{ vectors $u^{(j)}\in\R^m$ mutually orthogonal,}
\] 
\item $\alpha_1 v^{(1)} = \mA^\mT u^{(1)}$,
\[
 \mV_\ell = \left[\begin{array}{cccc} v^{(1)} & \cdots & v^{(\ell)}\end{array}\right], \quad \mbox{ vectors $v^{(j)}\in\R^n$ mutually orthogonal,}
\] 
\item $\mB_\ell\in\R^{(\ell+1)\times \ell}$ bidiagonal,
\[
 \mB_\ell = \left[\begin{array}{cccc} \alpha_1 &                 &          & \\
 \beta_2   & \alpha_2  &           & \\
  &\beta_3     &\ddots & \\
  &                 &\ddots & \alpha_\ell \\
  &                 &           &\beta_{\ell+1}
  \end{array}\right].
\]   
\item  The matrices satisfy the recurrence relations
\begin{eqnarray}
\mA \mV_\ell &=& \mU_{\ell+1} \mB_\ell, \label{AV} \\
\mA^\mT \mU_{\ell+1}&=& \mV_\ell\mB_\ell^\mT +\alpha_{\ell+1} v^{(\ell+1)} e_{\ell+1}^\mT, \label{ATU}
\end{eqnarray}                                                        
where $e_{\ell+1}\in\R^{\ell+1}$ is the canonical $(\ell+1)$th coordinate unit vector. 
\end{enumerate}

The algorithm in its basic form can be described in the following way.
\bigskip
\hrule
\medskip

{\bf Lanczos bidiagonalization algorithm}

\medskip
\hrule
\medskip

 \begin{enumerate}
 \item {\bf Given} $\mA \in\R^{m\times n}, \quad 0\neq b\in\R^m,\; 1<\ell <\min(m,n);$
 \item {\bf Initialize} $\beta_1= \|b\|;\; u^{(1)}= b/\beta_1; \;  \alpha_1 = \|\mA^\mT u^{(1)}\|;  \; v^{(1)} = \mA^\mT u^{(1)} /\alpha_1$.
 \item {\bf Repeat} for $j = 1,2,\ldots,\ell$
 \begin{enumerate}
 \item $\tilde{u}^{(j)} = \mA v^{(j-1)} - \alpha_{j-1}u^{(j-1)}, \quad \beta_j = \|\tilde{u}^{(j)}\|, \quad u^{(j)} = \tilde{u}^{(j)}/\beta_j$;
 \item $\tilde{v}^{(j)} = \mA^\mT u^{(j)} - \beta_j v^{(j-1)}, \quad \alpha_j = \|\tilde{v}^{(j)} \|, \quad v^{(j)} = \tilde{v}^{(j)}/\alpha_j$; 
 \end{enumerate}
 \item {\bf Compute}  $\tilde{u}^{(\ell+1)} = \mA v^{(\ell)} - \alpha_\ell u^{(\ell)},$ 
 \quad $\beta_{\ell+1}= \|\tilde{u}^{(\ell+1)}\|, $
 \quad $u^{(\ell+1)}= \tilde{u}^{(\ell+1)} /\beta_{\ell+1}$.
 \item {\bf Output:} Matrices  $\mU_{\ell+1}$, $\mB_\ell$, $\mV_\ell$.
 \end{enumerate}

\medskip
\hrule
\medskip

We observe that by construction, 
\[
 {\rm span}\{u^{(1)},\ldots,u^{(\ell)}\} = {\mathcal K}_\ell\big(b, \mA\mA^\mT \big), \quad 
  {\rm span}\{v^{(1)},\ldots,v^{(\ell)}\} = {\mathcal K}_\ell\big(\mA^\mT b, \mA^\mT\mA\big),
\]
where 
\[
{\mathcal K}_\ell\big(c, \mM \big)= {\rm span}\{c, \mM c, \ldots, \mM^{\ell-1}c\} 
\]
the $\ell$th Krylov subspace associated with the vector $c$ and the matrix $\mM$.  

The orthogonality of the two families of vectors Lanczos vectors $\{u^{(1)},\ldots,u^{(\ell)}\}$ and $\{v^{(1)},\ldots,v^{(\ell)}\}$ 
may not be retained in finite precision arithmetic as the number of Lanczos steps increases, due to the accumulation of roundoff errors. When the departure from orthogonality becomes significant, reorthogonalization may be used to restore it. 
Discussion of the numerical properties of the algorithm can be found in \cite{paige1982lsqr}. 

Next we will look at how the Lanczos process yields approximate solutions for the two linear systems of interest. Consider the adjoint problem (\ref{adjoint}) and an approximate solution in the form
\[
 z = \mU_{\ell} \xi, \quad \xi\in\R^{\ell},
\]
which can be written as
\[
 \mU_{\ell}\xi = \mU_{\ell+1}\widetilde\xi, \quad \widetilde\xi = \left[\begin{array}{c} \xi \\ 0\end{array}\right]\in \R^{\ell+1}.
\] 
Letting $\mC_\ell\in\R^{\ell\times\ell}$ be the matrix obtained by deleting the last row of the matrix $\mB_\ell$, it follows from (\ref{ATU}) that
\begin{eqnarray*}
 \mA^\mT z &=& \mA^\mT \mU_{\ell+1}\widetilde \xi = \mV_\ell\mB_\ell^\mT \widetilde\xi +\alpha_{\ell+1}v^{(\ell+1)}\widetilde \xi_{\ell+1} \\
 &=& \mV_\ell\mC_\ell^\mT\xi.
\end{eqnarray*}
Further, by using (\ref{AV}), we obtain
\[
 \mA\mA^\mT z =  \mA\mV_\ell\mC_\ell^\mT\xi  = \mU_{\ell+1}\mB_\ell\mC_\ell^\mT\xi,
\]   
where
\[
 \mB_\ell \mC_\ell^\mT\xi = \left[\begin{array}{c} \mC_\ell \\ \beta_{\ell+1}e_\ell^\mT\end{array}\right]\mC_\ell^\mT \xi =
 \left[\begin{array}{c} \mC_\ell\mC_\ell^\mT \xi \\ \beta_{\ell+1}e_\ell^\mT\mC_\ell^\mT \xi\end{array}\right].
\] 
Observe that
\[
 e_\ell^\mT\mC_\ell^\mT = \left[\begin{array}{cccc} 0 & \cdots & 0 &\alpha_\ell\end{array}\right],
\]
so 
\[
 \beta_{\ell+1}e_\ell^\mT\mC_\ell^\mT \xi = \beta_{\ell+1}\alpha_\ell\xi_\ell,
\] 
and we arrive at the conclusion that
\[
 \big(\mA\mA^\mT + \mu\,\mI_m) z = \mU_\ell\big(\mC_\ell\mC_\ell^\mT + \mu\,\mI_\ell\big)\xi +  \beta_{\ell+1}\alpha_\ell\xi_\ell u^{(\ell+1)}.
\] 
Moreover, since $b = \beta_1 u^{(1)}$, it follows from the orthogonality of the vectors $u^{(j)}$ that
\[
 \big\|\big(\mA\mA^\mT +\mu\, \mI_m\big)z - b\big\|^2 = \big\|\big(\mC_\ell\mC_\ell^\mT + \mu\,\mI_{\ell+1}\big)\xi - \beta_1 e_1\big\|^2 +
  \big(\beta_{\ell+1}\alpha_\ell\xi_\ell\big)^2,
\] 
thus the minimizer of this expression can be found by solving a linear least squares problem of size $(\ell+1)\times\ell$.
 
When solving the normal equations, we seek an approximate solution of the form
\[
 x_\ell = \mV_\ell \zeta,\quad \zeta\in\R^\ell,
\]
hence from (\ref{AV})
\[
 \mA x_\ell = \mA \mV_\ell \zeta =  \mU_{\ell+1}\mB_\ell \zeta
\]
and 
\[
 \mA^\mT\mA x_\ell = \mA^\mT\mU_{\ell+1} \mB_\ell\zeta  = \mV_\ell\mB_\ell^\mT\mB_\ell\zeta + \alpha_{\ell+1} v^{(\ell+1)} e_{\ell+1}^\mT\mB_\ell\zeta.   
\]
From the observation that 
\[
e_{\ell+1}^\mT\mB_\ell =  \left[\begin{array}{cccc} 0 & \cdots & 0 &\beta_{\ell+1}\end{array}\right],
\]
it follows that
\[
\big(\mA^\mT\mA +\mu\, \mI_n\big)x_\ell = \mV_\ell\big(\mB_\ell^\mT\mB_\ell +\mu\, \mI_\ell\big)\zeta +\alpha_{\ell+1}\beta_{\ell+1}\zeta_{\ell} v^{(\ell+1)}.
\]
Since the right hand side of the normal equations can be written as
\[
 \mA^\mT b = \beta_1\mA^\mT u^{(1)} =  \alpha_1\beta_1 v^{(1)},
\]
the orthonormality of the columns of $\mV_\ell$ implies that
\[
 \big\|\big(\mA^\mT\mA + \mu\,\mI_n\big)x - \mA^\mT b\|^2 =  \|\big(\mB_\ell^\mT\mB_\ell + \mu\, \mI_\ell\big)\zeta -\alpha_1\beta_1 e_1\big\|^2 +\big(\alpha_{\ell+1}\beta_{\ell +1}\zeta_{\ell}\big)^2,
\] 
thus the minimizer of the expression above can be found by solving an  $(\ell+1)\times \ell$ linear least squares problem.

We conclude by pointing out that the computation of the Lanczos iterates can carried out with minimal memory allocation, by discarding the previous orthonormal vectors and storing only the nonzero entries of the bidiagonal matrix. Once the small coefficients of the basis of the Krylov subspace have been determined, the basis vectors can be recomputed and added together on the fly.

\section{Priors and linear transformations}\label{sec:qr}

The discussion so far assumed a white noise model for the prior and the likelihood. More generally, let the distributions of $x$ and $\varepsilon$ be given by (\ref{gaussian model}) in Theorem~\ref{th:get xj},
and assume that the symmetric factorization of the precision matrices $\mGamma^{-1}$ and $\mSigma^{-1}$,
\[
 \mGamma^{-1} = \mL^\mT \mL, \quad \mSigma^{-1} = \mS^\mT \mS
\]
are available. After the change of variables
\[
 \widetilde x = \mL(x - x_0), \quad \widetilde\mA = \mS \mA \mL^{-1}, \quad b = \mS(b - \mA x_0),
\]
the problem can be cast in standard normal form (\ref{ATA perturbed}).
The impact of the whitening on the computational complexity depends on the properties of the matrices $\mL$ and $\mS$, and will be discussed in the context of the computed examples.

In computational inverse problems, it is fairly common to express a priori belief about the unknown solution in terms of its degree of smoothness.  Given a discretized vector $x$ of temporal or spatial values of a distributed parameter which is believed to be smooth, let $\mL$ be  discrete differential operator, e.g., based on finite difference approximation. Define the auxiliary variable  
\begin{equation}\label{difference}
 z = \mL x,
\end{equation}
and assume that $z\sim{\mathcal N}(0,\mI)$. Typically, $\mL\in\R^{p\times n}$, $p>n$, and $\mL$ is of full rank, implying that its null space is trivial, however whenever $\mL$ is not a square matrix, it cannot be regarded as a symmetric factor of a prior precision matrix. 

In general, smoothness priors can be in the form of a zero-mean white noise assumption for $z$,
\[
z\sim{\mathcal N}(0,\mI_p), \quad z \in{\mathcal R}(\mL),
\]
restricted to the range of $\mL$ to guarantee that $z$ is related to the variable $x$ of primary interest which can be retrieved by multiplying $z$ by the pseudoinverse of $\mL$, i.e.,
\[
 x = \mL^\dagger z.
\]
Given a linear observation model with additive Gaussian noise, after whitening of the noise, the posterior density is of the form
\[
 \pi_{z\mid b}(z\mid b) \propto {\rm exp}\left(-\frac12 \|b - \mA \mL^\dagger z\|^2 - \frac 12 \|z\|^2\right), \quad z \in {\mathcal R}(\mL).
\]
When sampling from this density, it is necessary to ensure that the sample 
$\{z^1,\ldots,z^N\}\subset\R^p$ is in the subspace 
${\mathcal R}(\mL)$, so that we can calculate the corresponding sample
\begin{equation}\label{x sample}
 \{x^1, \ldots, x^N\}\subset \R^n,\quad x^j = \mL^\dagger z^j.
\end{equation}
The sampling approach proposed above can be used in this context without significant modifications. In fact,  consider the optimization problem
\begin{equation}\label{minimizer z}
 z = {\rm argmin}\left\{\|b +\eta - \mA \mL^\dagger z\|^2 + \|z - \nu\|^2\right\}, \quad \eta \sim{\mathcal N}(0,\mI_m), \quad
 \nu \sim{\mathcal N}(0,\mI_p).
\end{equation}
After expressing both $z$ and $\nu$ as the sum of orthogonal components,
\[
 z= z_1 + z_2, \quad \nu = \nu_1+\nu_2, \quad z_1,\nu_1\in {\mathcal R}(\mL), \quad z_2,\nu_2\in {\mathcal N}(\mL^\mT),
\]
we can write the objective function to be minimized as
\[
\|b +\eta - \mA \mL^\dagger z\|^2 + \|z - \nu\|^2 =
\|b +\eta - \mA \mL^\dagger z_1\|^2 + \|z_1 - \nu_1\|^2 + \|z_2 - \nu_2\|^2.
\]
This implies that the minimizer is of the form
\[
 z = z_1+z_2 = z_1 + \nu_2, \quad z_1 = {\rm argmin}\left\{\|b +\eta - \mA \mL^\dagger z\|^2 + \|z - \nu_1\|^2\right\} \in{\mathcal R}(\mL),
\]
and, in view of the properties of the pseudoinverse,
\[
 x = \mL^\dagger z = \mL^\dagger z_1.
\]
In summary, the above reasoning shows that while the random draws $\nu^j \sim {\mathcal N}(0,\mI_p)$ typically contain a non-zero component orthogonal to the subspace ${\mathcal R}(\mL)$, the sampling scheme generates intermediate points $z^j\in\R^p$ with the same orthogonal component, which is projected to zero when the final sample (\ref{x sample}) is generated. Therefore, one can either project $\nu$ onto the subspace ${\mathcal R}(\mL)$ before computing $z_1$, or keep the orthogonal component in the intermediate sample and filter it out at the end by multiplying with the pseudoinverse of $\mL$, whichever is easier.

A final comment on the computation of the pseudoinverse is in order. In general, the pseudoinverse can be computed by the passing through the reduced QR factorization of $\mL$,
\[
 \mL = \mQ \mR, \quad \mQ\in\R^{p\times n}, \quad \mR \in\R^{n\times n},
\]
where $\mQ$ is a matrix with an orthogonal basis of ${\mathcal R}(\mL)$ in its columns and $\mR$ is invertible, and letting
\[
 \mL^\dagger = \mR^{-1}\mQ^\mT.
\]
Unfortunately, in high dimensions the QR decomposition may become a computational bottleneck. A computationally viable way of multiplying a vector by the pseudoinverse or its transpose is by observing that
\[
 \mL^\dagger z = \big(\mL^\mT\mL\big)^{-1}\mL^\mT z, \quad
 (\mL^\dagger)^\mT y = \mL \big(\mL^\mT\mL\big)^{-1} y.  
\]
While in general the sparsity of a matrix $\mL$ does not guarantee the sparsity of $\mL^\mT\mL$, in many applications $\mL$ is approximating local operators, e.g.,  differential operators, hence it is either sparse or have nonzero entries only on few diagonals and $\mL^\mT\mL$ is also sparse, as illustrated in the next section.

\section{Computed examples}

To demonstrate the proposed algorithm in practice, we consider four computed examples, two of which are versions of a linear inverse problems arising in tomography, followed by two non-linear inverse problems in which the Gaussian sampling is part of the work flow to generate proposals for an MCMC algorithm.

\begin{figure}[ht!]
\centerline{
\includegraphics[width=8cm]{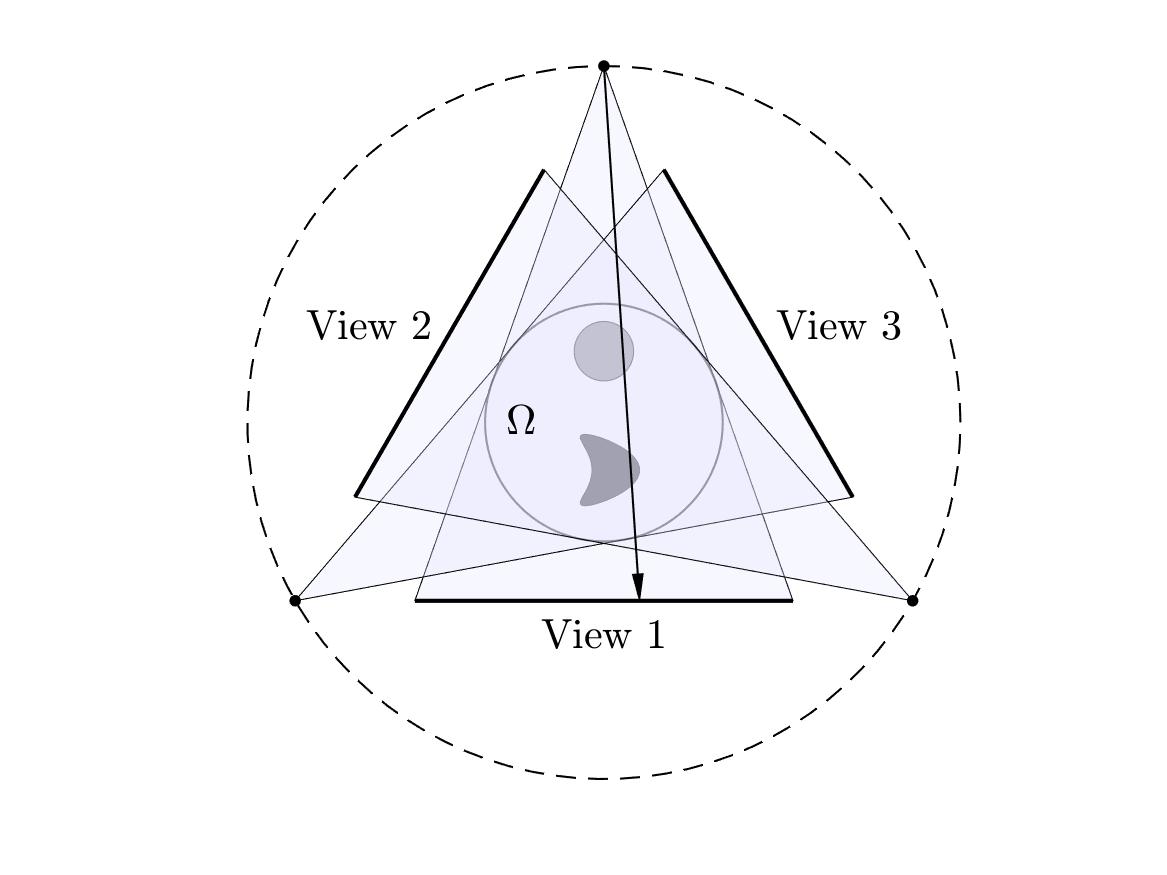}
}
\caption{The configuration of the fanbeam tomography example. In the figure, three equally spaced illumination views are shown.}\label{fig: tomo config}
\end{figure}

{\bf Fanbeam tomography with few projections.} In this example, we consider a fanbeam tomography problem schematically shown in Figure~\ref{fig: tomo config}. The circular object with unknown inclusions is illuminated from few directions by a point-like X-ray source and the transmitted projections are recorded at the opposite side by a linear array of receivers. The two-dimensional density distribution is reconstructed from the measured attenuation profiles along receiver arrays.

We assume that the data consist of noisy observations of the integrals of the density $\rho:\Omega\to\R$ over rays through the target, 
\begin{equation}\label{xray}
 b^*_j = \int_0^{\ell_j} \rho(\gamma_j(t)) dt,
\end{equation}
where $\gamma_j$ is the parametrization of the $j$th ray with respect to the arclength $t$ and $\ell_j$ is the length of the ray. We assume that the target is illuminated from 10 uniformly distributed directions, and each projection is recorded by 200 uniformly distributed receivers along the corresponding linear array, resulting in a data vector of dimension $m = 2\,000$.

\begin{figure}[ht!]
\centerline{\includegraphics[width=\textwidth]{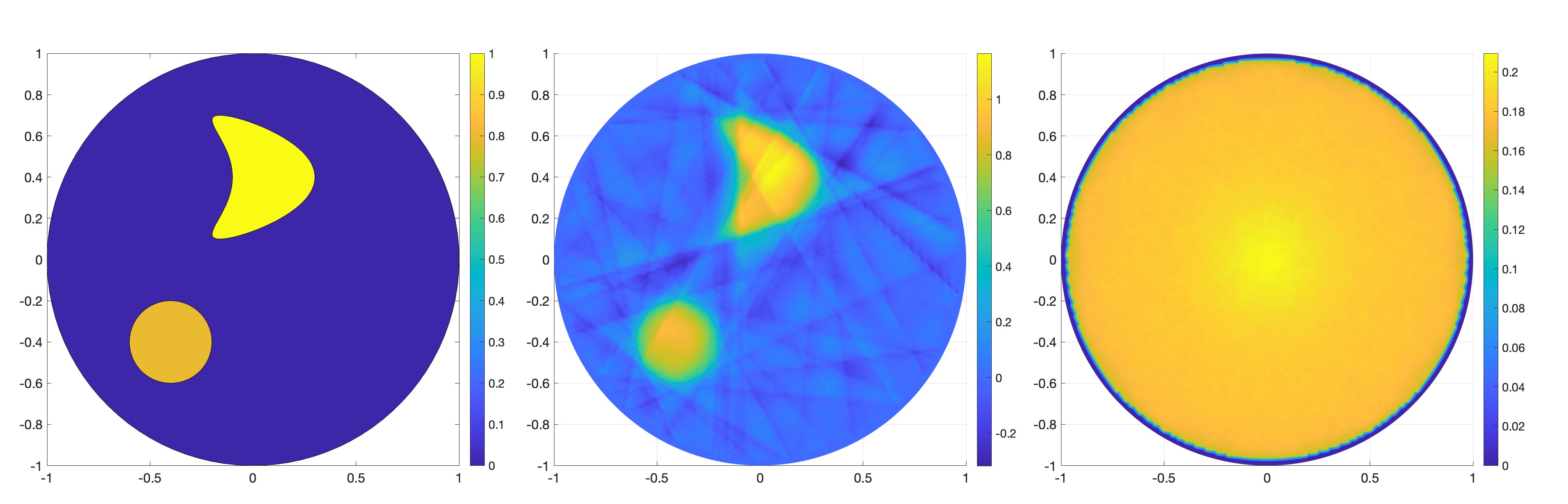}}
\caption{The generative model for the simulated data is shown in the left panel. The middle and the right panels show the estimated posterior mean and marginal standard deviations computed from a sample of size $10\,000$.}\label{fig: fanbeam}
\end{figure}

The data for our computed example are generated using a generative target shown in the left panel of Figure~\ref{fig: fanbeam} containing two inclusions. The data are generated in a mesh-free manner by computing the integrals along the rays, which can be done analytically as the boundary curves of the inclusions are analytically given. We assume that the background density of the target is zero, and that the constant densities in the circular and kite-shaped inclusions are 0.8 and 1, respectively.

For the inverse problem, we discretize the circular target, modeled as a unit disc, partitioning it into a triangular element mesh with $17\,813$ elements and $9\,062$ vertices, generated by using the DistMesh algorithm \cite{persson2004simple} with homogeneous mesh size parameter $h_{\rm mesh} = 0.02$. In the inverse problem. the density is parametrized by the nodal values over the mesh. We assume that the density of the target at boundary nodes is known, so the number of unknowns corresponding to the density values at the interior nodes is $n = 8\,455$. 

To define the prior, we introduce a first order finite difference matrix in the following way. Let $e_\ell = \{v_j,v_k\}$ denote an edge of the triangular mesh connecting the nodes $v_j$ and $v_k$, and let $p$ be the number of interior edges for which at least one endpoint is an interior node. We define the matrix $\mL\in\R^{p\times n}$ by setting
\[
 \mL_{\ell,j} =1,\quad \mL_{\ell,k} = -1,
\]
if both nodes are interior nodes. If one of the end nodes of the $\ell$th edge is a boundary node, the $\ell$th row of $\mL$ contains only one non-zero entry $\pm 1$ in the column corresponding to the node that is an interior node. In the current application, $p = 25\,640$, so the index of underdeterminacy is $m/p \approx 0.078$.
The matrix $\mL$ is scaled by a parameter $\alpha$ whose value is chosen experimentally so that the resulting penalty allows for jumps corresponding to targets with contrast of order of magnitude of unity. For this example, the value of the scaling parameter was set to $\alpha = 0.05$.
To formulate the problem in standard form, the multiplication by the pseudoinverse of $\mL$ was computed by using the formula 
\[
\mA \mL^\dagger = \mA (\mL^\mT \mL)^{-1} \mL^\mT = \big((\mL^\mT\mL)^{-1}\mA^\mT\big)^\mT\mL^\mT.
\]
Due to the sparsity of the matrix $\mL^\mT\mL$, this computation is very fast, requiring less than 4 seconds in a standard 2.3 GHz dual core laptop, or 0.7 seconds on a Mac Studio with Apple M2 Ultra chip.

The discretized forward model is assumed to be of the form (\ref{linsys}), where $\varepsilon$ denotes scaled white noise with standard deviation $\sigma^*$ equal to 10\% of the maximum noiseless data vector entry.  In the likelihood model, we used a value $\sigma$ slightly larger than the generative noise model to account for the discretization errors \cite{calvetti2020bayesian}.

To test the performance of the proposed algorithm, we generate samples of $100$, $500$, $5\,000$, and $10\,000$ draws using both the normal equation model (\ref{normal n}) and the proposed subspace splitting model, entailing the solution of the two adjoint equations (\ref{get delta}) and (\ref{adjoint m}).
The computing times using a direct solver on a Mac Studio platform with Apple M2 Ultra chip and 64 GB of memory are shown in Table~\ref{tab1}. Observe that since the sampling by direct methods can be done using multiple right hand sides with one matrix factorization, the computing times do not scale linearly as a function of the sample size, and the reduction in computing time using the latter method improves as the sample size grows. 

For comparison, we also computed the sample by solving the normal equations (\ref{normal n}) approximately using 30 steps of the Lanczos algorithm, yielding a relative accuracy of the solution in the range of $10^{-5} - 10^{-4}$. As expected, the computing times increase linearly with respect to the sample size, as the methodology does not solve the system with multiple right hand sides.

\begin{table}[ht!]
\begin{center}
\begin{tabular}{l|ccccc}
Sample size & 100 & 500 & 1000 & 5000 & 10000 \\
\hline
Normal eq., direct solver [s] & 124 & 167 & 225 & 681 &  1239\\
Adjoint eq., direct solver [s] & 19 & 21 & 22 & 40 & 60 \\
Ratio of adjoint/normal &15.3\% & 12.4 \% & 10.2\% & 6.2\% & 4.8\% \\
\hline
Lanczos solver & 89 & 439 & 877 & 4423 &  8819 \\
\end{tabular}
\end{center}
\caption{Computing times in seconds for different sample sizes using the algorithm based on the normal equations (\ref{normal n}), and on the proposed adjoint model with subspace splitting, including the solution of both (\ref{get delta}) and (\ref{adjoint m}). The ratio of the latter divided by the former is given in percentages. The last row lists the computing times using the iterative Lanczos solver for the normal equations. The platform is a Mac Studio with Apple M2 Ultra chip, 64 GB memory.}\label{tab1}
\end{table}

{\bf Cross-borehole tomography:} In the second example, we consider again a tomography problem but with scarcer data than in the previous model. Consider two parallel boreholes drilled into the Earth, and assume that a point-like seismic source is lowered into one of them, and the resulting seismic wave is recorded in the line of equally spaced receivers in the other borehole. It is assumed that the data at each receiver is a noisy observation of an integral of some physical property along the ray path between the source and the receiver. In travel time observations, the integrated parameter is the slowness of the wave, while the amplitude observation amounts to the integral of the attenuation, see, e.g., \cite{mcmechan1983seismic,ivansson1986seismic}. 
 
\begin{figure}
    \centerline{
\includegraphics[width=\textwidth]{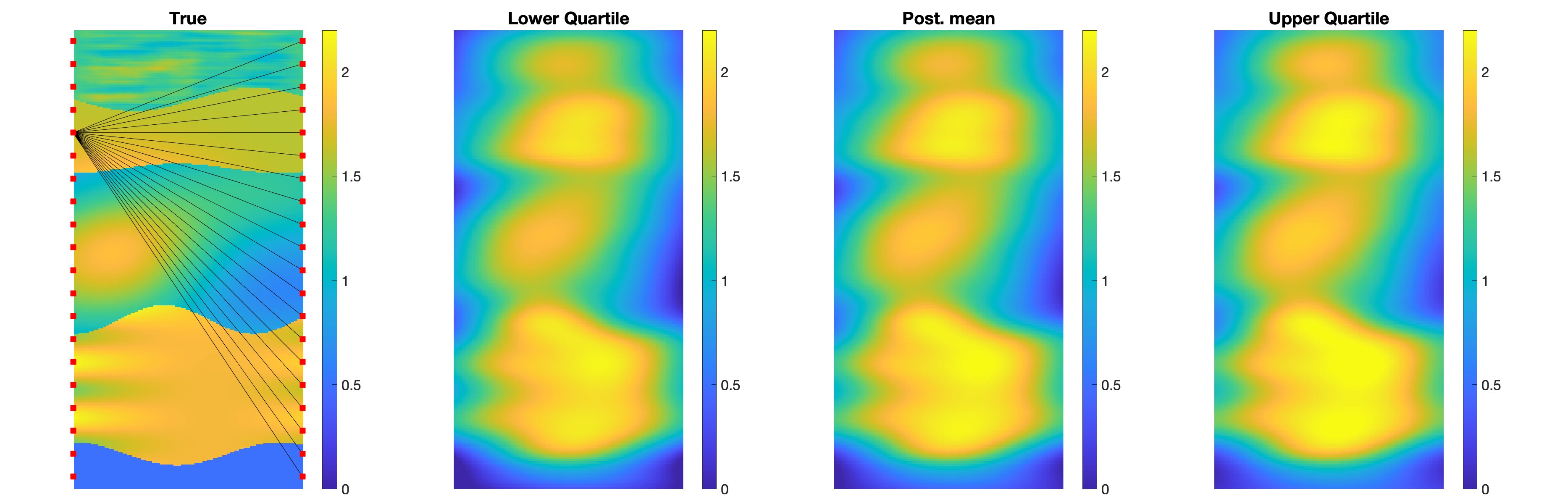}
}
\caption{Left: The generative structure with one seismic source in the left borehole, and a line of receivers in the right hole. The other three panels show the sample-based estimated posterior mean, and the lower and upper quartiles of 25\% and 75\%,  computed pixel-by-pixel.}\label{fig:borehole}
\end{figure}

The experiment is repeated with different source positions in the source borehole.  From the integral data along all possible ray paths, the pertinent property of the medium between the holes is retrieved. We refer to Figure~\ref{fig:borehole} for a graphical explanation of the setup.

To generate the data, we first define a structure consisting smooth variations and sharp almost horizontal jumps as shown in the left panel of Figure~\ref{fig:borehole}. We then discretize the domain by using a rectangular pixel mesh of size $300\times 140$ and approximate the structure by a pixel-wise constant function, allowing the integration along the rays. Gaussian additive scaled white noise is added to the data, with standard deviation of 0.5\% of the maximum entry of the noiseless data. For the solution of the inverse problem, we use a coarser model, dividing the rectangular domain of interest between the boreholes into $n= n_y\times n_z = 100\times 200 = 20\,000$ rectangular pixels and generate the tomography matrix assuming that the unknown structure is constant over pixels. The number of transmitter and receiver positions is assumed to be equal to 20 in each borehole, hence the dimension of the data space is $m = 20\times 20 = 400$.  In this problem, the index of underdeterminacy is $m/n = 400/20\, 000 = 0.02$.

To specify the prior model, we organize the pixel values in column-by-column order and use a Gaussian prior of anisotropic Whittle-Mat\'{e}rn type,
\[
  \mL X = W \sim {\mathcal N}(0,\gamma^2 \mI_n), 
\]
where $\mL = \mL_{\rm vert}\otimes \mL_{\rm hor}$ is the Kronecker product of the finite difference matrices in the vertical and horizontal directions, respectively,
\[
 \mL_{\rm vert} = \mL_{n_z} +\lambda_z^{-2}\mI_{n_z}, \quad \mL_{\rm hor} = \mL_{n_y} +\lambda_y^{-2}\mI_{n_y}.
\]
Here, $\mL_{n_z}\in\R^{n_z\times n_z}$ is a positive definite second order finite difference matrix in the vertical direction, and $\lambda_z$ is a vertical correlation length, and $\mL_{n_y}\in\R^{n_y\times n_y}$ is the corresponding matrix in the horizontal direction. We assume a priori that the unknown structure has longer correlation in the horizontal than in the vertical direction, and set $\lambda_x = 10$ and $\lambda_y = 20$, the units in pixels. Finally, $\gamma$ is a prior scaling adjusted so that the prior variance corresponds to the a priori expectations about the dynamical range of the image, $\gamma =70$. 

We run the proposed algorithm to generate a sample of size 10\,000 sample points, and for comparison, we compute a sample of the same size by using the algorithm based on the standard normal equations form. The corresponding computing times are $t_{\rm adj} = 5.5$ seconds and $t_{\rm normal} = 736$ seconds, yielding the ratio
\[
 \frac{t_{\rm adj}}{t_{\rm normal}}\times 100 \approx 0.7\%. 
\]

Figure~\ref{fig:borehole} shows the sample-based estimates for the posterior mean estimate, and the lower an upper pixel-wise quartiles. 
 
 {\bf Electrical impedance tomography:} In this example, we consider the non-linear electrical impedance tomography (EIT) inverse problem of estimating the conductivity distribution inside a body based on current-voltage measurements on its boundary. In this computed two-dimensional problem, $\Omega\subset\R^2$ is a bounded smooth domain with connected complement representing the body with unknown conductivity distribution, denoted by $\sigma$. We assume here that $\sigma$ is continuous up to the boundary $\partial\Omega$, and $\sigma(x)\geq \sigma_{\rm min}$ for some $\sigma_{\rm min}>0$. We assume that $L$ contact electrodes, modeled as non-overlapping closed intervals $e_\ell\subset\partial\Omega$ of positive length of the boundary curve, are attached on the boundary surface, and electric currents $I_\ell$ are injected through them. The current injections constitute a current pattern, a vector $I$ in $\R^L$. For each current pattern, the electric voltages over the electrodes are measured, giving rise to the corresponding voltage pattern $U\in \R^L$. The current patterns must satisfy the Kirchhoff's law, $\sum_{\ell=1}^L I_\ell=0$, thus we may choose $L-1$ linearly independent current patterns, $\{I^1,I^2,\ldots,I^{L-1}\}$, called a frame, and measure the corresponding voltages, $\{U^1,U^2,\ldots,U^{L-1}\}$. 

Given the conductivity distribution $\sigma$ and a current pattern $I$, the idealized forward model with the complete electrode model (CEM) consists of solving the elliptic PDE for the electric potential $u$ in $\Omega$,
\[
 \nabla \cdot \sigma\nabla u = 0,
\]
with the boundary conditions  
\[
 \left(u +z_\ell \sigma\frac{\partial u}{\partial n}\right)\bigg|_{e_\ell} = U_\ell, \quad \int_{e_\ell}\sigma\frac{\partial u}{\partial n} dS = I_\ell, \quad \frac{\partial u}{\partial n}\bigg|_{\partial\Omega \setminus \cap e_\ell} = 0,
\]  
where the constants $z_\ell>0$ are the contact impedances and presumably known. The above problem is a well-posed elliptic problem \cite{somersalo1992existence} that allows a stable approximation by the finite element method (FEM). Here, we assume that the domain $\Omega$ is approximated by a tessellation with triangular elements  $\{K_\nu\}$. We assume that in a thin collar domain near the boundary, the conductivity is constant $\sigma_0>0$ and known, while in the interior domain consisting of $n$ elements that we refer to as interior triangles, the conductivity may vary. Using a piecewise first order polynomial Lagrange basis, the numerical approximation depends on the conductivity only through the integrals of it over the interior triangles. We denote the  degrees of freedom by 
\[
 \xi_\nu = \frac 1{|K_\nu|}\int_{K_\nu} \sigma(x) dx, \quad 1\leq \nu\leq n.
\] 
To guarantee positivity of the conductivity, we reparametrize the conductivity and define new dimensionless variables $\gamma_\nu$ by setting
\[
  \xi_\nu = \sigma_0\big(e^{\gamma_\nu} - 1\big).
\]  
Given a fixed current frame and the tessellation of $\Omega$,  a discrete FEM approximation of the forward model is denoted by
\[
F:\R^n \to \R^m, \gamma \mapsto F(\gamma) =  \left[\begin{array}{c} U^1 \\ \vdots \\ U^{L-1}\end{array}\right], \quad m = L(L-1).
\]
For the numerical and computational details, we refer to \cite{calvetti2024sparsity}.  

For the computational test, we assume additive Gaussian noise model,
\[
 B = F(\Gamma) + E,
\]
where  $\Gamma$ is an $n$-variate random variable with realizations $\gamma$, and $E\sim{\mathcal N}(0, \eta^2 \mI_m)$, where $\eta>0$ is the noise level, leading to the likelihood model
\[
 \pi_{B\mid \Gamma}(b\mid\gamma)\propto {\rm exp}\left(-\frac 1{2\eta^2}\| F(\gamma) - b\|^2\right).
\] 
 
\begin{figure}[ht!]
\centerline{
\includegraphics[width=\textwidth]{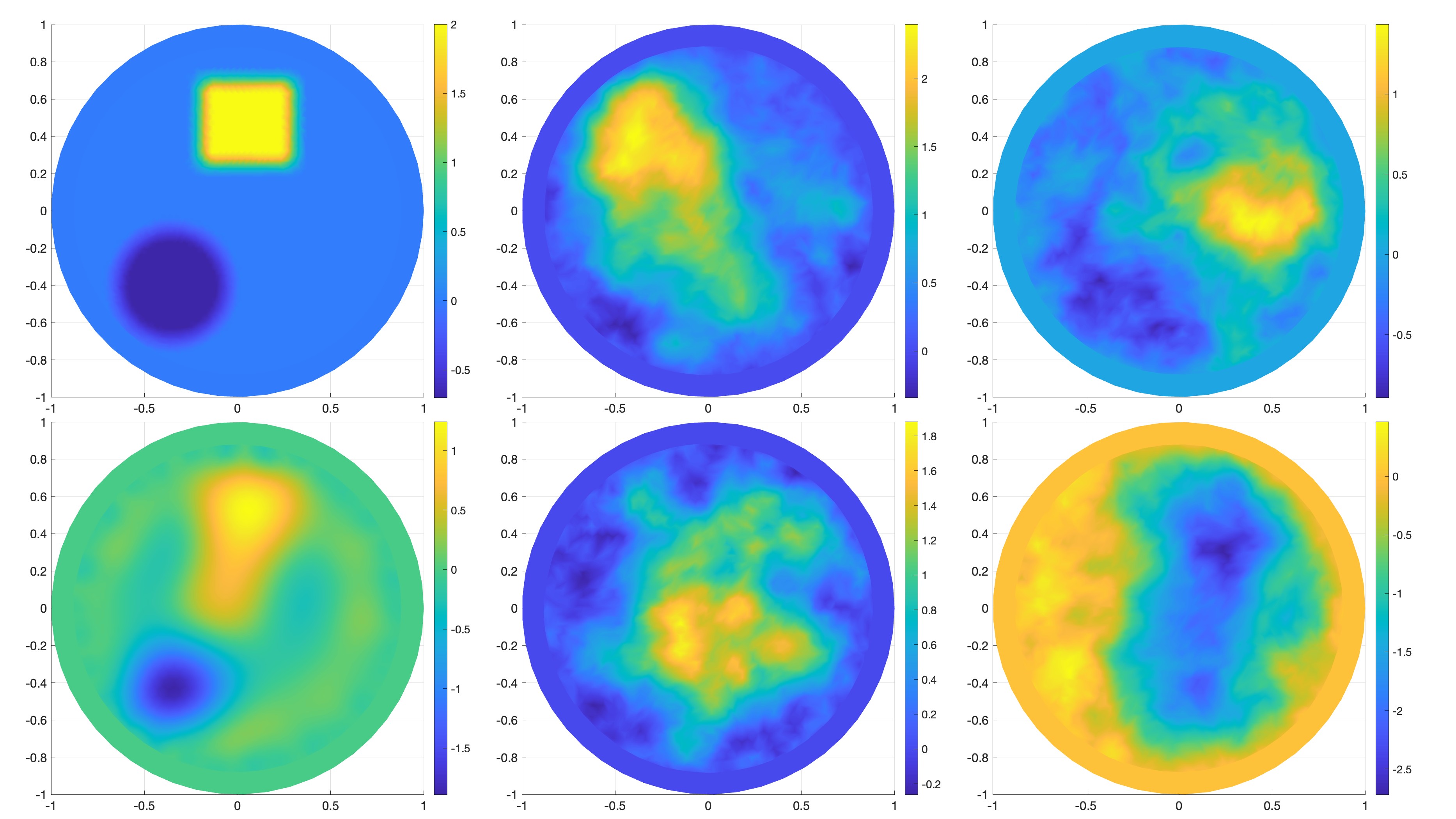}
}    
\caption{On the top left, the generative model $\gamma =\log(\sigma/\sigma_0+1)$ for the data is shown, and the bottom left panel shows the corresponding estimate with the linearized model. The four panels in the middle and on the right show four random draws from the prior distribution.}\label{fig:Generative and draws}
\end{figure}

For the prior, we assume that $\Gamma$ is slowly varying, described by a Whittle-Mat\'{e}rn model,
\[
 \pi_\Gamma(\gamma) \propto {\rm exp}\left( - \frac 12\| \mL \gamma\|^2\right),
\]
where $\mL\in\R^{n\times n}$ is a matrix of the form
\[
 \mL = \delta\left( \mD^\mT\mD  + \frac 1{\lambda^2} \mI_n\right),
\] 
where $\mD\in\R^{k\times n}$, $k>n$, is a first order finite difference matrix, $k$ being the number of edges of the mesh consisting of interior triangles. The mapping $\gamma\mapsto \mD\gamma$ returns the jumps across the edges, with the understanding that outside the triangle mesh, $\gamma$ vanishes. The parameter $\lambda>0$ is a correlation length parameter, and $\delta>0$ is a scaling parameter.  We set the  parameters $(\lambda,\delta)$ experimentally so that random draws from the prior reproduce details of the size that we a priori expect to see in the distribution $\gamma$, and that the marginal variances of the elements correspond to the expected dynamic range of the variable $\gamma$. To generate synthetic data that is in agreement with the prior model, we use a generative model shown in Figure~\ref{fig:Generative and draws}, where a few random draws from the prior are also shown.

As a first test, we consider the linearized model around the background conductivity $\sigma_0$,
\[
F(\gamma) \approx F(0) + DF(0) \gamma.
\]
By scaling the data and the forward map by $1/\eta$, we may write the likelihood approximately as
\[
 \pi_{B\mid\Gamma}(b\mid\gamma) \propto
{\rm exp}\left(-\frac 1{2}\| F(\gamma) - b\|^2\right) \approx {\rm exp}\left(-\frac 1{2}\| D F(0) \gamma - (b - F(0))\|^2\right)
\]
and, by introducing the notations $\mA = DF(0)$   and $r = b - F(0)$, we obtain the approximate Gaussian posterior
\[
 \pi_{\Gamma\mid B}(\gamma\mid b) \approx \pi^G_{\Gamma\mid B}(\gamma \mid b) \propto {\rm exp}\left( -\frac 12\|\mA\gamma - r\|^2 + \|\mL\gamma\|^2\right).
\]
The Gaussian sampler can be applied immediately to this distribution, leading to an independent sample from this approximate distribution. For later reference, we write the distribution as
\[
 \pi^G_{\Gamma\mid B}(\gamma\mid b) = {\mathcal N}(\gamma\mid \overline\gamma,\mC),
 \]
 where
 \[
 \overline\gamma = \left(\mA^\mT\mA + \mL^\mT\mL\right)^{-1}\mA^\mT b, \quad \mC =  \left(\mA^\mT\mA + \mL^\mT\mL\right)^{-1}.
 \]

To generate a sample from the fully non-linear posterior density, we write the linearization error as
\[
 M(\gamma) = F(\gamma) - F(0) - DF(0)\gamma,
\]
and the posterior density as
\[
 \pi_{\Gamma\mid B}(\gamma\mid b) \propto {\rm exp}\left(-\frac 1{2}\| \mA\gamma - r  +M\|^2\right) = \pi^G_{\Gamma\mid B}(\gamma\mid b)e^{-\Phi(\gamma)},
\]
where
\[
 \Phi(\gamma) = \frac12  \|M(\gamma)\|^2 + M(\gamma)^\mT(\mA\gamma - r).  
\] 
To generate a sample from the posterior using the machinery developed in this article, we employ the preconditioned Crank-Nicholson (pCN) sampler \cite{cotter2013mcmc} discussed below. We refer to
\cite{dunlop2015bayesian} for pCN applied to the EIT problem.

To initialize the sampler, we set $\gamma^0 = \overline\gamma$, and $\delta\gamma^0 = \gamma^0-\overline\gamma = 0$, and set the counter to $k=1$.

To generate sample $\gamma^k$ at the current value $k$, we draw a proposal move
\[
 w^* \sim {\mathcal N}(0, \mC),
\]
and given the previous sample point $\gamma^{k-1}$ with $\delta\gamma^{k-1} = \gamma^k - \overline\gamma$, we define the proposal
\[
 \gamma^* = \overline\gamma + \sqrt{1-h^2}\delta\gamma^{k-1} + h w^k,
\]  
where $0<h<1$ is a stepsize control parameter. To complete the step, we compute the acceptance ratio
\[
 \alpha = {\rm exp}\left( -\Phi(\gamma^*) + \Phi(\gamma^{k-1})\right),
\]
and accept the move with probability $\min(1,\alpha)$ setting $\gamma^k = \gamma^*$, otherwise rejecting it and setting $\gamma^k = \gamma^{k-1}$.

The philosophy of the sampler is based on the assumption that the linearization is already a relatively good approximation of the model, and the acceptance test depends mostly on the linearization error. Indeed, if we had $M(\gamma)=0$, the moves are automatically accepted.

Finally, to generate the proposal draws $w$, observe that by the form of the posterior density, we have
\[
 {\mathcal N}(w\mid 0,\mC) \propto {\rm exp}\left( -\frac 12\|\mA\gamma \|^2 -\frac 12 \|\mL\gamma\|^2\right),
\]
that is, we may use the same Gaussian sampler that generates the independent samples from $\pi^G_{\Gamma\mid B}(\gamma\mid b)$ by simply setting $r=0$. 

To test the sampler, we implemented a two-dimensional EIT model corresponding to $L=32$ equally spaced electrodes around a circular domain. The generative model for the synthetic data is shown in Figure~\ref{fig:Generative and draws}. In the model, we assume that in the annulus $\{0.9<|x|<1\}$, the conductivity has a known background value $\sigma=1$. To accommodate the electrode boundaries where the current densities have singularities, the triangular FEM meshes in the annular domain are made denser towards the boundary, while in the interior disc $V = \{|x|<0.9\}$, a homogeneous triangularization was used. The number of elements of the generative mesh was  10\,788, of which 5\,418 in the interior disc $V$, while in the mesh used for solving the inverse problem, the total number of elements was 5\,816, of which $n=1\,940$ are inside the domain $V$. The data were generated by using a full frame of trigonometric current feeds, the voltages being measured on every electrode, leading to the dimensionality of the data $m = L(L-1) = 992$, the index of underdeterminacy in this case being $m/n = 0.51$.

After forming the finite element matrices, we computed the Jacobian defining the matrix $\mA$, and generated 10\,000 random proposal draws for the pCN sampler. The computing time using the adjoint sampler on a Mac Studio was 0.9 seconds. We then run the pCN algorithm, using a stepsize $h = 0.05$. With this step size, the acceptance rate was 38\%. The most time consuming step in this process is the evaluation of the forward model, leading to a total computing time  of 22 minutes and 11 seconds.



{\bf Dictionary learning problem in magnetoencephalography:} In this problem, we consider the magnetoencephalography (MEG) data used to identify an active region of a brain atlas based on non-invasive magnetic field measurement outside the head. Brain activity is characterized by neuronal firing activity: when tens of thousands of bundled neurons are simultaneously activated, a magnetic field detectable outside the skull is induced by the post-synaptic ion currents. The MEG inverse problem is to locate the activity from the field measurements. We refer to \cite{hamalainen1993magnetoencephalography,baillet2001electromagnetic} for details. The standard computational model approximates the neuronal activity by a sum of current dipoles representing the spatially localized post-synaptic ion currents of the neuron bundles, and the magnetic field as a sum of the magnetic fields induced by each dipole. Due to the linearity of Maxwell's equations, the vector of magnetic flux density values measured by $m$ magnetometers, denoted by $b$, is linearly related to the dipole moments $q_j$ of the $n$ dipoles, leading to a model
\[
 b = \sum_{j=1}^n\mL_j q_j, \quad q_j\in\R^3,
\]
where $\mL_j$ are the $ m \times 3$ lead field matrices that take into account the geometry and conductivity structure of the head \cite{sarvas1987basic,de2012eeg}. Typically, the number of dipoles is tens of thousands, while the number of measured magnetic field components is of the order of few hundreds.
An often used simplification assumes that the dipoles have a preferential direction, typically normal to the cortex, i.e., $q_j = x_j \gamma_j$, where $x_j\in\R$ is the unknown amplitude and $\gamma_j\in\R^3$ is a unit vector. With this simplification, the model can be written as
\[
 b = \sum_{j=1}^n \big(\mL_j\gamma_j\big) x_j = \mL x, \quad x \in\R^n, \quad \mL \in\R^{m\times n}.
\]
The MEG data consist of noisy observations of the vector $b$.

When studying the time-dependent brain activity, e.g., in order to map the brain connectivity during a resting state, it is often sufficient to identify the active functional brain region. Assuming that the source space of $n$ dipoles is subdivided into $L$ parcels, each parcel corresponding to $n_\ell$ dipoles with $n_1 + \ldots + n_L = n$, we may write the forward model as
\[
 b = \sum_{\ell = 1}^L \mL^\ell x^\ell,
\]
where $\mL^\ell\in\R^{m\times n_\ell}$, $x^\ell\in \R^{n_\ell}$. To identify the active brain region, we wish to find a solution such that most of the vectors $x^\ell$ vanish. This objective leads naturally to a Bayesian prior model promoting group sparsity. More precisely, we introduce a conditionally Gaussian prior density,
\[
 x^\ell\mid \theta_\ell \sim {\mathcal N}(0,\theta_\ell \mI_{n_\ell}), \quad \theta_\ell \sim{\rm InverseGamma}(\beta,\vartheta_\ell),
\]
where $\beta$ and $\vartheta_\ell$ are the shape and scaling parameters for the variance $\theta_\ell$. By assuming additive Gaussian scaled white noise as a measurement noise, the model leads to a posterior density of the form
\begin{eqnarray*}
&&\pi_{X,\Theta\mid B}(x,\theta\mid b) \\
&&\propto{\rm exp}\left( -\frac 1{2\sigma^2}\|b - \sum_{\ell = 1}^L \mL^\ell x^\ell\|^2 -\frac 12
\sum_{\ell=1}^L \frac{\|x^\ell\|^2}{\theta_\ell} - \sum_{\ell=1}^L \frac{\vartheta_\ell}{\theta_\ell} -\sum_{\ell=1}^L\left(\beta+1 + \frac{n_\ell}{2}\right)\log\theta_\ell\right),
\end{eqnarray*}
A viable approach for identifying the active brain region is to search for the maximum a posteriori (MAP) estimate, which identifies the brain regions that best explain the data \cite{pragliola2022overcomplete,waniorek2023bayesian}.
A computationally efficient algorithm for the MAP estimate,
the iterative alternating sequential (IAS) algorithm, proceeds as follows:
\begin{enumerate}
\item Given the parameters $\sigma$, $\vartheta_\ell$ and $\beta$, and the data $b$,
\item {Initialize} $\theta=\vartheta$,
\item {Iterate} until convergence:
\begin{enumerate}
\item  Update $x$ by solving the least squares problem
\[
 x = {\rm argmin} \left\{\frac 1{2\sigma^2}\|b - \sum_{\ell = 1}^L \mL^\ell x^\ell\|^2 +\frac 12
\sum_{\ell=1}^L \frac{\|x^\ell\|^2}{\theta_\ell}\right\}.
\]
\item Update $\theta$ from the component-wise the optimality conditions
\[
 \frac{\partial}{\partial \theta_\ell} \left\{
\frac 12
 \frac{\|x^\ell\|^2}{\theta_\ell} +  \frac{\vartheta_\ell}{\theta_\ell} +\left(\beta+1 + \frac{n_\ell}{2}\right)\log\theta_\ell
 \right\} = 0,
\]
yielding
\[
 \theta_\ell = \frac{1}{\kappa_\ell}\left(\frac{\|x^\ell\|^2}{2} + \vartheta_\ell\right), \quad \kappa_\ell = \beta+1 + \frac{n_\ell}{2}.
\]
\end{enumerate}
\end{enumerate}

The convergence properties of the algorithm have been investigated in a series of articles \cite{calvetti2019hierachical,calvetti2020sparse,pragliola2022overcomplete}. With the selection of an inverse Gamma hyperprior for the parameter $\theta$, the algorithm strongly promotes sparsity,  however, since the corresponding objective function is non-convex, there is no guarantee of a unique minimizer.
To keep the algorithm from stopping too early at a non-representative local minimizer, a hybrid version of the algorithm was proposed in \cite{calvetti2020sparsity}. A relevant question from the Bayesian viewpoint is how strongly the hierarchical prior is concentrated around sparse solutions, and how representative of the posterior distribution the MAP estimate is. One way to answer these questions is to explore the posterior density using Markov chain Monte Carlo (MCMC) sampling, see \cite{calvetti2024computationally}. In this application, we consider a block-Gibbs sampler, leveraging the proposed sampler from a Gaussian density as a subroutine of the sampler. More specifically, w propose the following sequential sampling algorithm.

\begin{enumerate}
\item Given: the parameters $\sigma$, $\vartheta_\ell$ and $\beta$, the data $b$, and $T$ = desired size of the sample.
\item Initialize $\theta = \vartheta$. Set the counter $t=1$. 
\item Generate a new sample point $(x,\theta)^t$, $1\leq t\leq T$ as follows:
\begin{enumerate}
\item Draw $x$ from the conditional posterior distribution
\[
 \pi_{X\mid B,\Theta}(x\mid b,\theta) \propto
 {\rm exp}\left( -\frac 1{2\sigma^2}\|b - \sum_{\ell = 1}^L \mL^\ell x^\ell\|^2 -\frac 12
\sum_{\ell=1}^L \frac{\|x^\ell\|^2}{\theta_\ell} \right),
\]
\item For $\ell = 1,\ldots,L$, draw $\theta_\ell$ from
\[
 \theta_\ell \sim{\rm InverseGamma}\left(\vartheta_\ell+\frac{\|x^\ell\|^2}{2},\beta + \frac{n_\ell}{2}\right).
\]
\end{enumerate}
\item If $t<T$, advance the counter $t\to t+1$, else, stop.
\end{enumerate}
The implementation of step 3 (a) is effectuated by using the proposed Gaussian sampler.

In our example, the discretized head model consist of $n = 15\, 002$ dipoles, and $m=306$ measurements performed by 102 sensors each comprising a magnetometer and a gradiometer, recording the magnetic field density as well as tangential gradient components. The head model and corresponding lead field matrix  $\mL \in\R^{m\times n}$ was downloaded from the Open MEG Archive (OMEGA) website {\tt https://www.mcgill.ca/bic/neuroinformatics/omega},  see \cite{niso2016omega}. The source space is divided into $L = 148$ parcels according to the Destrieux brain atlas \cite{destrieux2010automatic}. To test the algorithm, we randomly select one brain region, then randomly pick a  dipole in it. Starting at the selected dipole, we generate an activity patch by picking its four nearest neighbors within the same brain region, and we assign a dipole moment randomly drawn from a uniform distribution to each of the five dipoles. The data thus generated is then corrupted by adding scaled Gaussian white noise with standard deviation 0.5\% of the maximal component of the noiseless data. Next, we compute the MAP estimate using the IAS algorithm, and run the block-Gibbs sampler as described above. Observe that the Gibbs sampler is independent of the MAP estimator algorithm, i.e., we do not initialize it at the MAP but rather at the point $\theta = \vartheta$.
We refer to \cite{Calvetti2025dictionary} for details concerning the parameter selection as well as the discussion of the IAS algorithm, limiting our attention here to the computational performance of the sampling algorithm. 

\begin{figure}[ht!]
\includegraphics[width=\textwidth]{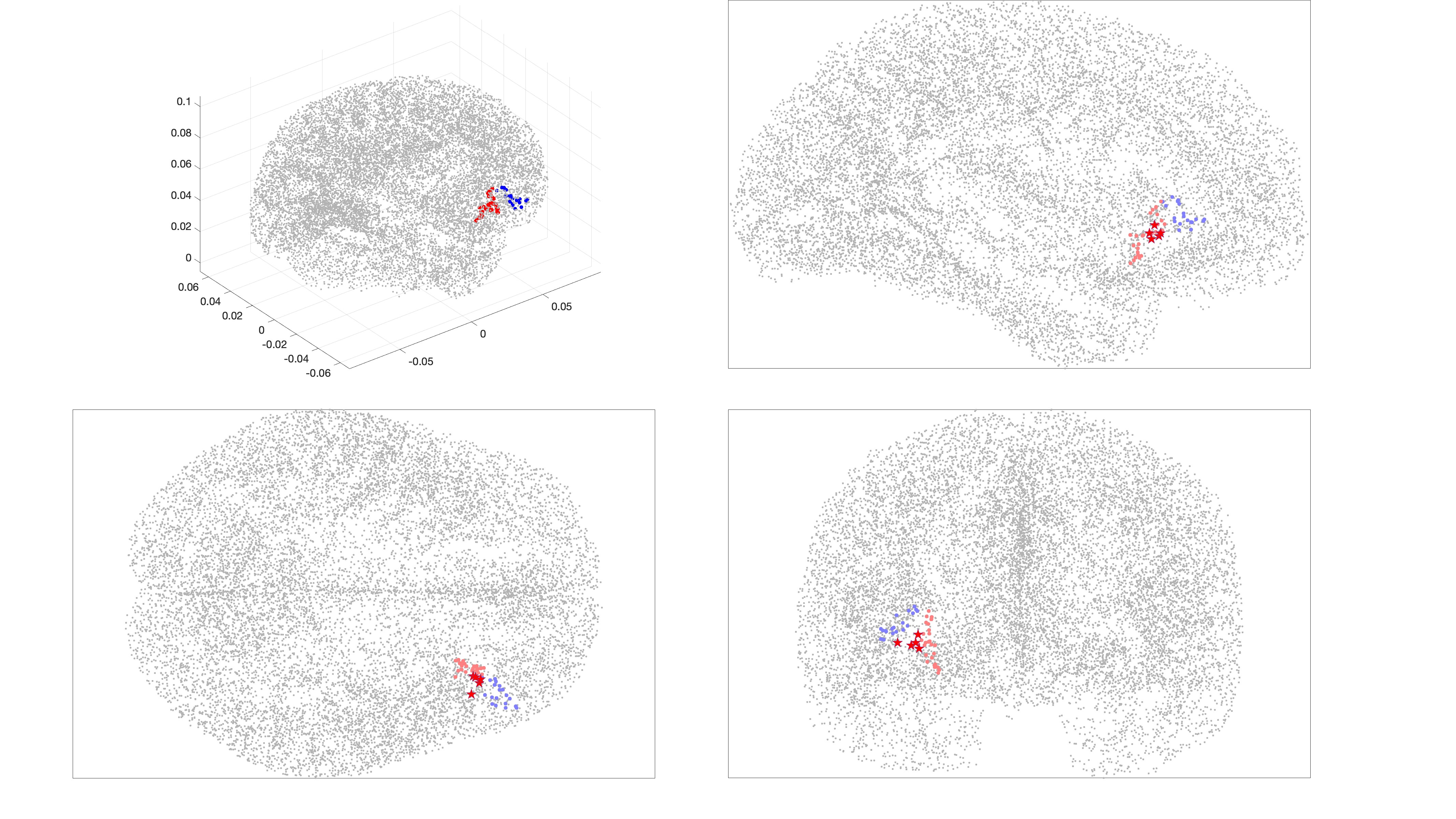}
\caption{The source space consisting of 15\,002 dipoles, their locations marked by gray dots. The selected active parcel of dipoles are located on the anterior part of the right circular insular sulcus, are indicated in red, while the active dipoles generating the data with stars. The blue dots indicate the dipole positions on the anterior horizontal part of the lateral fissure, corresponding to the parcel of dipoles having the second largest posterior expectation for the variance.}
\label{fig:brain1}
\end{figure}

\begin{figure}[ht!]
\includegraphics[width=\textwidth]{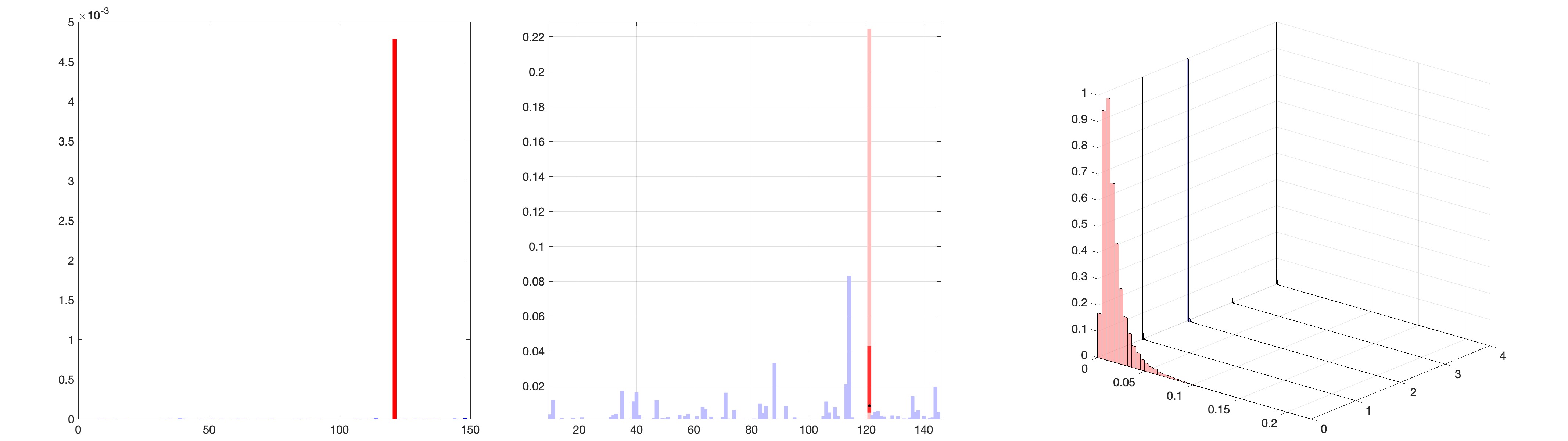}
\caption{Left panel shows the MAP estimate of the vector $\theta$ obtained by the IAS algorithm. The red bar corresponds to the correct brain region that the data was generated with. The middle panel shows the ranges of the sample values $\theta_j$ (light red for the correct brain region, light blue for the incorrect ones), together with the 75\% quantiles indicated by the dark red/dark blue). The right panel shows the histograms of the variances of the five brain regions corresponding to the largest posterior means of the variance parameters. The red histogram corresponds to the correctly identified brain region.}\label{fig:hist1}
\end{figure}

In the first numerical test, we generate activity in the anterior part of the right circular insular sulcus, which is a relatively deep source region. The source dipoles in this region are indicated by red dots in Figure~\ref{fig:brain1}, the active dipoles marked with stars. The IAS algorithm identifies the active region without ambiguity, as seen in Figure~\ref{fig:hist1}, left panel.  We then generate a sample of size 30\,000  drawn from the posterior density. The middle panel shows a quantile plot of the samples of the variance parameters, and the right panel shows the histograms of the five brain regions with the largest posterior mean variances. The second largest variance correspond to anterior horizontal  part of the lateral fissure, which is indicated by blue dots in Figure~\ref{fig:brain1}, confirming that the confounding domains are geometrically close to the true source. The histograms in Figure~\ref{fig:hist1} confirm that the IAS-MAP estimate is representative of the posterior, and rather reliable. Most importantly for the message of this article, the computing time for the sampling result is 18 minutes 45 seconds on a Mac Studio desktop computer. Observe that since the Gibbs sampler requires the solution in every iteration with a different matrix, one cannot take advantage of matrix factorizations with multiple right hand sides, and the computing time scales linearly as a function of the sample size.  

\section{Discussion}

In this work, we propose a computationally efficient sampling method for generating draws from a Gaussian posterior distribution when the data space has a significantly lower dimensionality than the space of the unknown. While the idea of solving normal equations related to the posterior mean estimation of linear inverse problems is well-known, in random sampling the equivalence of the normal equations and the adjoint equations breaks down, and an extra adjoint problem needs to be solved to obtain a subspace splitting that allows the use of the equivalence. 
Computed examples demonstrate the efficiency of the algorithm when the dimension of data space is small enough to allow the use of direct solvers. Computed examples show that the computational gain increases faster than the sample size, which is understandable as the direct methods can take advantage of solving the problem with multiple right hand sides with a single matrix factorization. When the matrix is not available or is too large to be stored in memory, thus requiring the use of iterative linear solvers, to there is no advantage in resorting to the adjoint problem, because the same Lanczos process is required to solve both the normal equations and the adjoint problem. A detailed computational complexity analysis addressing also memory considerations will be the topic of future work.

\section*{Acknowledgements}

The work of DC was partly supported by the NSF grant DMS 1951446, and that of  ES by the NSF grant DMS-2204618.

{\bf Conflict of interest:} The authors have no competing interests to declare that are relevant to the content of this article.

{\bf Data availability:} No datasets were generated or analyzed during the current study.

\bibliographystyle{siam} 
\bibliography{biblio}
\end{document}